\documentclass[a4paper,12pt]{amsart}
\usepackage[all]{xy}
\usepackage{amsmath, amsfonts, amssymb}
\usepackage{graphics}

\newtheorem{proposition}{Proposition}[section]
\newtheorem{theorem}[proposition]{Theorem}

\newtheorem{lemma}[proposition]{Lemma}

\newtheorem{definition}{Definition}[section]

\newcommand{\pc}{\ensuremath{\mathbb{Z}/p\mathbb{Z}}}

\newcommand{\s}[1]{\ensuremath{S^{#1}}}

\newcommand{\ad}{\ensuremath{\mathrm{ad}}}
\newcommand{\lra}{\longrightarrow}
\newcommand{\dr}[3]{\ensuremath{#1\stackrel{#2}
{\longrightarrow}#3}}
\newcommand{\ddr}[5]{\ensuremath{#1\stackrel{#2}
{\longrightarrow}#3\stackrel{#4}{\longrightarrow}#5}}
\newcommand{\dddr}[7]{\ensuremath{#1\stackrel{#2}
{\longrightarrow}#3\stackrel{#4}{\longrightarrow}#5
\stackrel{#6}{\longrightarrow}#7}}
\newcommand{\ddddr}[9]{\ensuremath{#1\stackrel{#2}
{\longrightarrow}#3\stackrel{#4}{\longrightarrow}#5
\stackrel{#6}{\longrightarrow}#7}\stackrel{#8}{\longrightarrow}#9}

\newcommand{\hahc}{homotopy associative, homotopy commutative }

\renewcommand{\L}{\ensuremath{\mathcal{L}}}

\begin{document}
\title{Universal spaces of two-cell complexes and their exponent bounds}
\author{Jelena Grbi\' c}
\maketitle
\begin{abstract}
Let $P^{2n+1}$ be a two-cell complex which is formed by attaching a
$(2n+1)$--cell to a $2m$--sphere by a suspension map. We construct a
universal space $U$ for $P^{2n+1}$ in the category of homotopy
associative, homotopy commutative $H$--spaces. By universal we mean
that $U$ is homotopy associative, homotopy commutative, and has the
property that any map $f\colon P^{2n+1}\lra Y$ to a homotopy
associative, homotopy commutative $H$--space $Y$ extends to a
uniquely determined $H$--map $\overline{f}\colon U\lra Y$. We then
prove upper and lower bounds of the $H$--homotopy exponent of $U$.
In the case of a mod~$p^r$ Moore space $U$ is the homotopy fibre
$S^{2n+1}\{p^r\}$ of the $p^r$--power map on $S^{2n+1}$, and we
reproduce Neisendorfer's result that $S^{2n+1}\{p^r\}$ is homotopy
associative, homotopy commutative and that the $p^r$--power map on
$S^{2n+1}\{p^r\}$ is null homotopic. \smallskip
\end{abstract}

\section{Introduction}
\label{intro}

Free or universal objects have been of interest to mathematicians in
many mathematical disciplines. In homotopy theory one of the first
universal spaces is given in the category of homotopy associative
$H$--spaces by the James construction~(cf.~\cite{Jr}). The analogue
of the James construction on a topological connected space $X$ in
the category of \hahc $p$--localised $H$--spaces is given by the
following definition.
\begin{definition}
\label{defuni} Localised at $p$ a \textit{universal space $U_X$} of
a topological space $X$ is a \hahc $H$--space together with a map
$i\colon X\lra U_X$ such that the following \textit{property} holds:
\begin{itemize}
\item [] if $Y$ is a \hahc localised at $p$ ${H\text{--space}}$
and ${f\colon X\lra Y}$ is any map, then $f$ extends to a unique
${H\text{--map}}$ ${\overline{f}\colon U_X\lra Y}$.
%As with the James construction, $\overline{f}$
%is called \textit{the multiplicative extension of} f.
\end{itemize}
\end{definition}
This whole concept of studying universal spaces in the category of
homotopy associative, homotopy commutative $H$--spaces is due to
Gray~(cf.~\cite{GEHP}).

Despite the potential applications of universal spaces $U_X$ in
homotopy theory, there are as yet just a few known examples of these
spaces. As listed in \cite{JGu}, universal spaces are known to exist
in the following cases. For one-cell complexes, spheres; for
two-cell complexes, Moore spaces and the $(2np-2)$--skeleton $K$ of
$\Omega^2S^{2n+1}$; and for a three-cell complex L which is
$(2np-1)$--skeleton of $\Omega^2S^{2n+1}$ are known.

The objective of this paper is to add to the list of examples a
family of universal spaces of certain two-cell complexes. We
construct these universal spaces using the method developed
in~\cite{JGu} and rely heavily on Gray's decomposition of loop
spaces on certain two-cell complexes~\cite{Go}. The basic underlying
reason why the methods given in~\cite{Go}, \cite{JGu} work is the
existence of a differential Lie algebra structure on homotopy groups
with coefficients.

Our standing assumptions are that all spaces have a non-degenerate
basepoint, are simply connected, have the homotopy type of a
$CW$--complex, and are localised at an odd prime $p>3$. Unless
otherwise indicated, the ring of homology coefficients will be $\pc$
and $H_*(X;\pc )$ will be written as $H_*(X)$.

The main object of consideration is given as follows.
\begin{definition}
For $n\geq m$, let $ \Theta\colon S^{2n-1}\lra S^{2m-1}$ be a given
map. Then the two-cell complex $P$ is defined as the homotopy
cofibre of~$\Theta$.
\end{definition}
We decorate $P$ with a superscript to denote the dimension of the
top cell. Namely, $P^r=\Sigma^{r-2n}P$ for $r\geq 2n$.

Let $U$ be the homotopy fibre of $\dr{\Sigma^2\Theta\colon
S^{2n+1}}{}{S^{2m+1}}$. Then our main theorem is:
\begin{theorem}
\label{uniP}
 $U$ is a universal space of $P^{2n+1}$.
\end{theorem}

The uniqueness assertion of an $H$--extension in Theorem~\ref{uniP}
is powerful. It ensures that if a universal space exists then it and
its $H$--structure that is homotopy associative and homotopy
commutative are unique up to homotopy equivalence. The uniqueness of
an $H$--extension can be also used to show that two $H$--maps $U\lra
Y$ are homotopic by comparing their restrictions on $P^{2n+1}$.
Further more it is important to point out that a two-cell complex
$P^{2n+1}$ is given as a mapping cone of an attaching map of finite
order. We use this to obtain exponent information about the spaces
$U$ in Theorem~\ref{uniP}.

For an arbitrary space $Y$ we define its \emph{homotopy exponent},
denoted $\exp(Y)=p^t$, if $t$ is the minimal power of $p$ which
annihilates the $p$--torsion in the homotopy groups of $Y$.

\begin{theorem}
\label{lubounds} Let $p^l$ be the order of $\Sigma\Theta$ and $p^k$
the order of $\Sigma^2\Theta$. Then
\[p^{n-k}\leq \exp(U)\leq p^{n+l}.\]
\end{theorem}

As a special case of our work we reprove some results of
Neisendorfer (cf.~\cite{Nprop}) related to the mod~$p^r$ Moore space
$P^{2n+1}(p^r)$. If the degree $p^r$ map on $S^{2n}$ is taken for
$\Sigma\Theta$, then the corresponding two-cell complex is mod~$p^r$
Moore space $P^{2n+1}(p^r)$. According to Theorem~\ref{uniP}, a
universal space of $P^{2n+1}(p^r)$ is the homotopy fibre of
$\dr{\Sigma^2\Theta\colon S^{2n+1}}{}{S^{2n+1}}$. As the degree
$p^r$ map and the $p^r$--power map on an odd dimensional sphere
coincide, a universal space of $P^{2n+1}(p^r)$ is $S^{2n+1}\{p^r\}$,
the homotopy fibre of the $p^r$--power map on $S^{2n+1}$. Thus
$S^{2n+1}\{p^r\}$ is a \hahc $H$--space.

The mod--$p$ $H$--\emph{exponent} of an $H$--space $X$, denoted
$H\exp(X)$, is $p^k$ if $k$ is the minimal power for which the
$p^{\text{th}}$--power map on $X$ is null homotopic. We use the
universal property of $S^{2n+1}\{p^r\}$ to find the $H$-- exponent
of $S^{2n+1}\{p^r\}$.

\begin{proposition}
$H\exp(S^{2n+1}\{p^r\})=p^r$.
\end{proposition}
Neisendorfer's proofs of $H$--structure and $H$--exponent of
$S^{2n+1}\{p^r\}$ involved turning homotopy fibrations into actual
fibrations, and so were point set in nature. Using universal spaces,
our proofs retain the flexibility of the homotopies and so are
perhaps more transparent.
\bigskip

\noindent\textsc{Acknowledgments.} The author would like to thank
Brayton Gray for directing her to his work.
\section{Preliminaries}
In this section we recall some definitions and known results that
will be used in subsequent sections.

Define the map $\dr{\ad\colon\bigvee_{i=0}^\infty
P^{4n+3+2mi}}{}{P^{2n+2}}$ as the wedge sum of generalised Whitehead
products given by
$$
\ad=\bigvee^\infty_{i=0}\ad_i:P^{4n+3}\wedge \s{2mi}\lra P^{2n+2}
$$
where $\ad_i=[\iota , \ad_{i-1}]$ and
$\ad_0=[1_{P^{2n+2}},1_{P^{2n+2}}]\circ\Sigma q$, while $q$ is a
right homotopy inverse to ${P^{2n+1}\wedge P^{2n+1}\lra P^{4n+2}}$
and ${\iota\colon S^{2m+1}\lra P^{2n+2}}$ is the inclusion of the
bottom cell in $P^{2n+2}$. In \cite[Theorem 1.2]{Go} Gray described
the homotopy fibre of the map ad and gave a decomposition of $\Omega
P^{2n+2}$. Recall that $U$ is the homotopy fibre of
${\dr{\Sigma^2\Theta\colon S^{2n+1}}{}{S^{2m+1}}}$.
\begin{theorem}
\label{fibrationP} There exists a homotopy fibration sequence:
\begin{equation}
\label{fibP} \ddddr{\Omega \Big(\bigvee_{i=0}^\infty
P^{4n+3+2mi}\Big)}{\Omega \ad}{\Omega
P^{2n+2}}{\partial}{U}{*}{\bigvee_{i=0}^\infty
P^{4n+3+2mi}}{\ad}{P^{2n+2}}
\end{equation}
and a homotopy decomposition
\begin{equation}\label{gspliting}
\Omega P^{2n+2}\simeq U\times\Omega\Big(\bigvee_{i=0}^\infty
P^{4n+3+2mi}\Big).
\end{equation}
\end{theorem}
\subsection{Homotopy Groups with $P$--coefficients}

Homotopy groups with coefficients in $P$ are defined as
\[
\pi_k(X;P)=[P^k, X] \mbox{  for } k\geq 0.
\]
Gray showed (cf.~{\cite [Proposition 3.5]{Gas}}) that there is a Lie
algebra structure on the homotopy groups with the coefficients in
$P$ given by mod~$P$ Samelson products. To begin, there is a
splitting
\begin{equation}
\label{splitP}
 P^k\wedge P^l \simeq  P^{k+l-2n+2m-1}\vee P^{k+l}.
\end{equation}
Let $\mu^*_{k,l}\colon P^{k+l}\lra P^k\wedge P^l$ be the inclusion.

If $G$ is a group-like space, and $\alpha\in \pi_k(G;P),$
${\beta\in\pi_l(G;P)}$, then the \emph{mod~$P$ Samelson product}
$[\alpha,\ \beta]$ is defined as the composition
$$
\ddr{P^{k+l}}{\mu^*_{k,l}}{P^k\wedge P^l}{\langle\alpha
,\beta\rangle}{G}
$$
where $\langle \alpha ,\beta\rangle$ is the commutator
$\langle\alpha,\beta\rangle=\alpha\beta\alpha^{-1}\beta^{-1}$ in
$G$. \emph{Mod~$P$ Whitehead products} are defined as the adjoints
of mod~$P$ Samelson products. Gray showed that $\mu^*_{k,l}$ can be
chosen so that the mod $P$ Samelson product satisfies anti-symmetry
and Jacobi identities.

Let $\sigma =2n-2m+1$ be the dimensional gap between the top and the
bottom cell of the two-cell complex $P$. Gray defined a Bockstein
homomorphism
\[
\beta\colon \pi_k(X;P)\lra\pi_{k-\sigma}(X;P)
\]
for $k\geq4n-2m+1$ that is a degree $\sigma$ derivation with respect
to the Lie bracket on the homotopy groups with coefficients in $P$.
If $X$ is a group-like space, $u\in\pi_k(X;P)$ and
$v\in\pi_l(X;P)$,~then
\[
\beta [u,v]=[\beta u,v]+(-1)^k[u,\beta v].
\]

The reduced homology $\widetilde{H}_*(P^k;\pc)$ is a free
$\pc$--module on two generators $v$ and $u$ in degrees $k$ and
$k-\sigma$, respectively. The Hurewicz map
\[
h\colon\pi_k (X;P)\lra \widetilde{H}_k(X;\mathbb{Z}_{(p)}).
\]
is given by $h(\alpha )=f_*(v)$, where $f$ is in $\pi_k(X;P)$ and
$f\colon P^k\lra X$ represents~$\alpha$. It is a Lie algebra
homomorphism (cf.~\cite{Gas}) from the homotopy groups with $P$
coefficients of $X$ to the Lie algebra of primitives in
$H_*(X;\pc)$, that is, if $X$ is a group-like space and
${\alpha\in\pi_*(X;P)}$, $\beta\in\pi_*(X,P)$, then
\begin{equation}
\label{HurSamel} h([\alpha ,\beta])=[h(\alpha ),h(\beta )].
\end{equation}

\subsection{A homology decomposition of $\Omega\Sigma P^{2n+1}$}

Notice that $H_*(P^{2n+1})$ is the free $\pc$--module with basis
$x,u$ in degrees $2m$ and $2n+1$. Since $H_*(P^{2n+1})$ is a trivial
coalgebra, $H_*(\Omega\Sigma P^{2n+1})$ is the primitively generated
tensor Hopf algebra $T(x, u)$ generated by $x$ and $u$. Therefore
$H_*(\Omega\Sigma P^{2n+1})$ can be considered as the universal
enveloping Hopf algebra $U\L$ of the free Lie algebra $\L=\L\langle
x, u\rangle$. Denote the commutator of $\mathcal{L}$ by
$[\mathcal{L} ,\mathcal{L}]$ and regard $\L_{ab}\langle x,u\rangle$
as the free graded abelian Lie algebra, that is, as $\L/ [\L,\L]$.
The geometrical decomposition~\eqref{gspliting} implies there is an
isomorphism in homology
\begin{equation}
\label{Liedecomposition} U\mathcal{L}\cong
U[\mathcal{L},\mathcal{L}]\otimes U\L_{ab}\langle x, u \rangle
\end{equation}
of left $U[\mathcal{L}, \mathcal{L}]$--modules and right
$U\L_{ab}\langle x,u \rangle$--comodules. A basis for $[\mathcal{L},
\mathcal{L}]$ is given in~\cite{Go} by
\begin{equation}
\label{basis} W=\{ y_k=\ad^k(x)[u,u]\text{ and }z_k=\ad^k(x)[x,u]\}.
\end{equation}
Let $\widetilde{W}$ be the set of maps $\ad_k\colon P^{4n+3+2mk}\lra
P^{2n+2}$ and their Bocksteins. Pair of basis elements $(y_k,z_k)$
are the mod~$P$ Hurewicz images of the maps $\ad_k$ and their
Bocksteins respectively. The $H$-map
$\dr{\Omega(\bigvee_{i=0}^\infty P^{4n+3+2mi})}{\Omega\ad}{\Omega
P^{2n+2}}$ is uniquely determined by its restriction to
$\bigvee_{i=0}^\infty P^{4n+2+2mi}$, and so the image of
$(\Omega\ad)_*$ is $U[\mathcal{L}, \mathcal{L}]$.

\section{Homotopy Associativity and Homotopy Commutativity}

To begin the argument proving the homotopy associativity and
homotopy commutativity of $U$ we recall a pair of constructions from
classical homotopy theory.

Let $X$ be a topological space, and $i_L$ and $i_R$ the inclusions
of $X$ into the wedge ${X\vee X}$. Looping, we can take the Samelson
product $[\Omega i_L, \Omega i_R]$. Adjointing gives the Whitehead
product $[\zeta_L,\zeta_R]$, where ${\zeta_A=i_{A}\circ ev}$ for
$A=L,R$ and $ev$ is the canonical evaluation map $ev\colon \Sigma
\Omega X\longrightarrow X$. A classical result in homotopy theory
asserts that there is a homotopy fibration
\begin{equation}
\label{classical} \ddr{\Sigma\Omega X \wedge \Omega
X}{[\zeta_L,\zeta_R]}{X\vee X}{i}{X\times X}
\end{equation}
where $i$ is the inclusion. The \textit{universal Whitehead product
of} $X$ is defined as the composition
\[
\Phi\colon\Sigma\Omega X \wedge \Omega
X\stackrel{[\zeta_L,\zeta_R]}{\lra} X\vee X\stackrel{\nabla}{\lra} X
\]
where $\nabla$ is the fold map. Notice that any Whitehead product on
$X$ factors through the universal Whitehead product of $X$.

The \textit{universal Samelson product of} $X$ is defined as the
adjoint of the universal Whitehead product of $X$, namely, as the
commutator of the identity map on $\Omega X$
\[
[i_{\Omega X},i_{\Omega X}]\colon \Omega X \wedge\Omega X\lra \Omega
X.
\]
The linchpin in showing that $U$ is \hahc is the following theorem
proved by Theriault (cf.~\cite{Tp}).
\begin{theorem}
\label{htpy assoc&commut} Let $\Omega
B\stackrel{\partial}{\longrightarrow} F \longrightarrow
E\longrightarrow B$ be a homotopy fibration sequence in which
$\partial$ has a right homotopy inverse. Suppose that there is a
homotopy commutative diagram
\[
\xymatrix{\Sigma\Omega B\wedge\Omega B \ar[r] \ar[d] & B\vee B
\ar[d]^{\nabla} \\
E \ar[r] & B, }
\]
where the upper composite in the square is the universal Whitehead
product of $B$. Then the multiplication on $F$ induced by the
retraction off $\Omega B$ is both \hahc and the connecting map
$\partial$ is an $H$--map.
\end{theorem}

Returning to our case, consider the fibration sequence
\[
\dddr{\Omega P^{2n+2}}{\partial}{U}{*} {\bigvee_{i=0}^\infty
P^{4n+3+2mi}}{\ad}{P^{2n+2}}
\]
of Theorem \ref{fibrationP}. Let $\dr{\Psi\colon\Sigma\Omega
P^{2n+2}\wedge\Omega P^{2n+2}}{}{P^{2n+2}\vee P^{2n+2}}$ be the
Whitehead product $[\zeta_L,\zeta_R]$ for $X=P$.
\begin{lemma}
\label{p-sumP} The Whitehead product
$$
\dr{\Psi\colon\Sigma\Omega P^{2n+2}\wedge\Omega
P^{2n+2}}{}{P^{2n+2}\vee P^{2n+2}}
$$
is homotopic to a sum of mod~P Whitehead products.
\end{lemma}
\begin{proof}
Using James' theorem and splitting \eqref{splitP}, $\Sigma\Omega
P^{2n+2}\wedge\Omega P^{2n+2} \simeq \Sigma M$ for $M$ a wedge of
two-cell complexes $P^k$. Let $\theta$ be the restriction of
$\Omega\Psi$ to $M$
\[
\dddr{\theta\colon M}{E}{\Omega\Sigma M}{\simeq}{\Omega(\Sigma\Omega
P^{2n+2}\wedge\Omega P^{2n+2})}{\Omega\Psi}{\Omega(P^{2n+2}\vee
P^{2n+2}).}
\]
Then $\theta$ is a Samelson product as it is the adjoint of the
Whitehead product $\Psi$. Since $\Omega\Psi$ has a left homotopy
inverse, as fibration \eqref{classical} splits when looped, it is an
inclusion in homology. Furthermore, the Hurewicz image of each
summand $P^k$ of $M$ under the composite $\theta$ is a bracket in
$\L\langle V\rangle$, where $H_{\ast}(\Omega (P^{2n+2}\vee
P^{2n+2}))\cong U\L\langle V\rangle$ and
$V=\widetilde{H}_\ast(P^{2n+1}\vee P^{2n+1})$. Using the identity
and mod~$P$ Bockstein maps on each summand of $P^{2n+2}\vee
P^{2n+2}$, it is clear that there exists a mod~$P$ Samelson product
on $\Omega(P^{2n+2}\vee P^{2n+2})$ which has the same Hurewicz image
as $P^k$. Summing these mod~$P$ Samelson products, one for each
summand of $M$, gives a map ${\lambda\colon M\lra\Omega(P^{2n+2}\vee
P^{2n+2})}$ with the property that $\lambda_*=\theta_*$. Each
mod~$P$ Samelson product factors through the loop of the universal
Whitehead product of $P$, as every Whitehead products on $P$ factors
through the universal Whitehead product of $P$. Therefore $\lambda$
lifts to a map $\lambda'\colon M\lra \Omega(\Sigma\Omega
P^{2n+2}\wedge \Omega P^{2n+2})$ with ${\lambda\simeq
\Omega\Psi\circ\lambda'}$. Extend $\lambda'$ to
$\dr{\overline{\lambda}\colon\Omega\Sigma M}{}{\Omega(\Sigma\Omega
P\wedge\Omega P)}$. As $\lambda_*=\theta_*$, we get
$(\Omega\Psi\circ\overline{\lambda})_*=(\Omega\Psi)_*$. As
$(\Omega\Psi)_*$ is a monomorphism, we must have
$(\overline{\lambda})_*$ is an isomorphism. Hence
$\overline{\lambda}$ is a homotopy equivalence. Taking adjoints then
proves the Lemma.
\end{proof}
By definition, the universal Whitehead product of $P^{2n+2}$ is the
composition $\ddr{\Phi\colon\Sigma\Omega P^{2n+2}\wedge\Omega
P^{2n+2}}{\Psi}{P^{2n+2}\vee P^{2n+2}}{\nabla}{P^{2n+2}}$.
\begin{lemma}
\label{liftP}
 There is a lift
\[
\xymatrix{ & \Sigma\Omega P^{2n+2}\wedge\Omega P^{2n+2}
\ar@{-->}[dl]\ar[d]^{\Phi}\\
\bigvee_{i=0}^\infty P^{4n+3+2mi} \ar[r]^-{\ad} & P^{2n+2}}
\]
of the universal Whitehead product $\Phi$ of $P^{2n+2}$ to
$\bigvee_{i=0}^\infty P^{4n+3+2mi}$.
\end{lemma}
\begin{proof}
Lemma \ref{p-sumP} shows that the universal Whitehead product on
$P^{2n+2}$ is homotopic to a sum of mod~$P$ Whitehead products. The
mod~$P$ Whitehead product defines the Lie bracket on the homotopy
groups with $P$ coefficients, endowing it with a Lie algebra
structure. The set $\widetilde{W}$, defined after~\eqref{basis},
consists of mod~$P$ Whitehead products which form a Lie basis for
$[\L,\L]$. So the mod~$P$ Whitehead products from the universal
Whitehead product can be rewritten as a linear combination of basis
elements.
\end{proof}

\begin{theorem}
\label{Fshahc} $U$ is a homotopy associative, homotopy commutative
$H$--space.
\end{theorem}
\begin{proof}
Applying Theorem \ref{htpy assoc&commut} to the fibration sequence
$$
\dddr{\Omega P^{2n+2}}{\partial}{U}{*} {P^{4n+3}\rtimes\Omega
S^{2m+1}}{\ad}{P^{2n+2}}
$$
and using Lemma \ref{liftP}, the Proposition follows.
\end{proof}

\section{A Universal Property of $U$}

In this section we show that $U$ satisfies the universal property in
the category of homotopy associative, homotopy commutative
$H$--spaces. Let $f\colon P^{2n+1}\lra Z$ be a map into a homotopy
associative, homotopy commutative $H$--space. We show that there is
a unique multiplicative extension $\overline{f}\colon U\lra Z$ of
$f$. The following Proposition is due to Gray (cf.\cite{GEHP}),
although the proof itself is adjust to the notation used in this
paper.

\begin{proposition} \label{trivialmapP} Let
$h\colon\Omega P^{2n+2} \lra Z$ be an $H$--map into a homotopy
commutative $H$--space $Z$. Then it factors through
$\partial\colon\Omega P^{2n+2}\lra U$.
\end{proposition}

\begin{proof}
Let $g$ be a right homotopy inverse of the $H$--map $\partial$ and
let $e$ be the homotopy equivalence ${e\colon \dr{U\times
\Omega\Sigma R}{g\cdot\Omega \ad}{\Omega P^{2n+2}}}$ (cf.
Theorem~\ref{fibrationP}). Define two maps $a,b\colon\Omega
P^{2n+2}\lra \Omega P^{2n+2}$ by the composites
\[
\dddr{a\colon\Omega\Sigma P^{2n+1}}{e^{-1}}{U\times \Omega\Sigma
R}{\pi_1}{U}{g}{\Omega\Sigma P^{2n+1}}
\]
\[
\dddr{b\colon\Omega\Sigma P^{2n+1}}{e^{-1}}{U\times \Omega\Sigma
R}{\pi_2}{\Omega\Sigma R}{\Omega \ad}{\Omega\Sigma P^{2n+1}}.
\]
In the following diagram
\[
\xymatrix{\Omega P^{2n+2}\ar[r]^-\Delta \ar[rd]^{e^{-1}}& \Omega
P^{2n+2}\times \Omega P^{2n+2}\ar[r]^-{a\times b} & \Omega P^{2n+2}
\times\Omega P^{2n+2} \ar[r]^-{\mu} & \Omega P^{2n+2}\\
& U\times \Omega\Sigma R \ar[r]^-{g\times\Omega \ad} & \Omega
P^{2n+2}\times\Omega P^{2n+2}\ar@{=}[u] &}
\]
the composition along the top row is $a+b$, while the bottom row is
the identity map on $\Omega P^{2n+2}$. The commutativity of the
diagram gives $\mathrm{Id}_{\Omega P^{2n+2}}\simeq a+b$.

Being an $H$--map, $h$ is determined by its restrictions on each of
the factors of $\Omega P^{2n+2}$, that is, $h\simeq h\circ
\mathrm{Id}_{\Omega P^{2n+2}}\simeq h\circ(a+b)\simeq h\circ a
+h\circ b$.

If the composite
\[
\ddr{\bigvee_{i=0}^\infty P^{4n+3+2mi}}{\Omega \ad}{\Omega
P^{2n+2}}{h}{Z}
\]
is null homotopic, then we have $h\circ b\simeq *$ and hence
$h\simeq h\circ a$ proving the Proposition. As $h\circ\Omega \ad$ is
the composite of $H$--maps, it is itself an $H$--map. Therefore by
the James construction, it is uniquely determined by its restriction
to $\bigvee_{i=0}^\infty P^{4n+2+2mi}$. The composite
$$
\ddr{\bigvee_{i=0}^\infty
P^{4n+2+2mi}}{E}{\Omega\Sigma\big(\bigvee_{i=0}^\infty
P^{4n+2+2mi}\big)}{\Omega \ad}{\Omega P^{2n+2}}
$$
is a wedge of mod~$P$ Samelson products as it is the adjoint of the
wedge of mod~$P$ Whitehead products
$\dr{\ad\colon\Sigma\big(\bigvee_{i=0}^\infty
P^{4n+2+2mi}\big)}{}{P^{2n+2}}$. Being an $H$--map, $h$ preserves
Samelson products. Therefore the wedge of mod~$P$ Samelson products
$(\Omega \ad)\circ E$, composed with $h$ into the homotopy
commutative $H$--space $Z$ is trivial.
\end{proof}

\begin{theorem}
\label{universalFs} Let $Z$ be a \hahc $H$--space. Let $f\colon
P^{2n+1}\lra Z$ be given. Then $f$ extends to an $H$--map
$\overline{f}\colon U\lra Z$, which is unique up to homotopy.
\end{theorem}
\begin{proof}
Consider the fibration sequence $\ddr{\Omega
P^{2n+2}}{\partial}{U}{*}{\bigvee_{i=0}^\infty P^{4n+3+2mi}}$.
Define the map $\overline{f}$ as the composite
\[
\ddr{\overline{f}\colon U}{g}{\Omega P^{2n+2}}{\tilde{f}}{Z},
\]
where $g\colon U\lra \Omega P^{2n+2}$ is a right homotopy inverse of
${\partial\colon\Omega P^{2n+2}\lra U}$ and $\tilde{f}\colon\Omega
P^{2n+2}\lra Z$ is the canonical multiplicative extension of $f$
given by the James construction. Look at the map $\overline{f}$ as a
candidate for the multiplicative extension of $f$. There is a
commutative diagram
\begin{equation}
\label{extension} \xymatrix{ P^{2n+1} \ar[r]^-E\ar[dr]^f &
\Omega\Sigma P^{2n+1} \ar[d]^{\widetilde{f}}\ar[r]^-{\partial} & U
\ar[dl]^{\overline{f}}\\
& Z &}
\end{equation}
where the left triangle commutes by the James construction, while
the commutativity of the right triangle is given by
Proposition~\ref{trivialmapP}. Diagram~\eqref{extension} ensures
that $\overline{f}$ is an extension of $f$.

Now we shall prove that $\overline{f}$ is an $H$--map by showing
that the diagram
\begin{equation}
\label{HmapP} \xymatrix{U \times U \ar[r]^-{g\times g}
\ar[dr]_-{\overline{f}\times\overline{f}} & \Omega P^{2n+2}
\times\Omega P^{2n+2}
\ar[r]^-{\mu}\ar[d]^-{\tilde{f}\times\tilde{f}} & \Omega
P^{2n+2}\ar[r]^-{\partial}\ar[d]^-{\tilde{f}} &
U\ar[d]^-{\overline{f}}\\
& Z\times Z \ar[r]^-{\mu} & Z\ar@{=}[r] & Z }
\end{equation}
commutes. The left triangle commutes by definition; the middle
square commutes since $\tilde{f}$ is an $H$--map; and the
commutativity of the right square is given by
Proposition~\ref{trivialmapP}. Summing this up,
diagram~\eqref{HmapP} commutes.

Finally we are left to show that $\overline{f}=\tilde{f}\circ
g\colon U\longrightarrow Z$ is the unique $H$--map extending
$f\colon P^{2n+1}\longrightarrow Z$. To prove this we use the
uniqueness of $\tilde{f}$ asserted by the James construction and the
result of Theorem \ref{htpy assoc&commut} which establishes that the
fibration connecting map ${\partial\colon\dr{\Omega P^{2n+2}}{}{U}}$
is an $H$--map. Let $\overline{f}$, $\overline{l}$ be two extensions
$\dr{U}{\overline{f},\ \overline{l}}{Z}$ of the map $f\colon\dr
{P^{2n+1}}{}{Z}$ which are $H$--maps. Precompose both maps with the
$H$--map $\partial\colon\dr{\Omega P^{2n+2}}{}{U}$. We obtain two
multiplicative extensions $\ddr{\Omega P^{2n+2}}{\partial}{U}
{\overline{f},\ \overline{l}}{Z}$ of $f\colon\dr{P^{2n+1}}{}{Z}$. By
the uniqueness of an $H$--map $\Omega P^{2n+2}\lra Z$ extending $f$,
it follows that $\overline{f}\circ\partial\simeq \overline{l}
\circ\partial$. Precomposing both compositions with the right
homotopy inverse $g\colon\dr{U}{}{\Omega P^{2n+2}}$ of the map
$\partial$, we get
$$
\overline{f}\circ\partial\circ g\simeq
\overline{l}\circ\partial\circ g.
$$
Hence
$$
\overline{f}\simeq \overline{l}
$$
and the uniqueness assertion is proved. This finishes the proof of
the Theorem.
\end{proof}
The Theorems~\ref{Fshahc} and~\ref{universalFs} together prove our
main result, Theorem~\ref{uniP}.

\section{The exponent for $U$}
\subsection{The $H$--exponent for $S^{2n+1}\{p^r\}$}

We first reprove Neisendorfer's result (cf.~\cite{Nprop}) that
$H\exp(S^{2n+1}\{p^r\})=p^r$ by using the universality of
$S^{2n+1}\{p^r\}$.
\begin{theorem}
$H\exp(S^{2n+1}\{p^r\})=p^r$
\end{theorem}
\begin{proof}
Recall that the degree $p^r$ map on the mod~$p^r$ Moore space is
null homotopic. Consider the map
\[
f\colon \ddr{P^{2n+1}(p^r)}{p^l}{P^{2n+1}(p^r)}{i}{S^{2n+1}\{p^r\}}
\]
given as the composition of the degree $p^l$ map on $P^{2n+1}(p^r)$
with the inclusion of the bottom two cell into $S^{2n+1}\{p^r\}$. As
$f$ is a map into a homotopy associative, homotopy commutative
$H$--space, Theorem~\ref{uniP} says that it can be extended to a
unique $H$--map $\overline{f}\colon U\lra S^{2n+1}\{p^r\}$ where $U$
is a universal space of $P^{2n+1}(p^r)$. But $U\simeq
S^{2n+1}\{p^r\}$. Any $k$--power map on $S^{2n+1}\{p^r\}$ is given
by the composite
$\ddr{S^{2n+1}\{p^r\}}{\Delta_k}{\big(S^{2n+1}\{p^r\}\big)^k}{\mu}{S^{2n+1}\{p^r\}}$.
Since $S^{2n+1}\{p^r\}$ is homotopy commutative, the multiplication
$\mu$ is an $H$--map and so it is the $k$--power map on
$S^{2n+1}\{p^r\}$. Taking the $p^r$--power map on $S^{2n+1}\{p^r\}$
for $\overline{f}$, there is a homotopy commutative diagram
\[
\xymatrix{ P^{2n+1}(p^r)\ar[d]\ar[r]^{{p^l}} & P^{2n+1}(p^r)\ar[r]^{i} & S^{2n+1}\{p^r\} \\
S^{2n+1}\{p^r\} \ar@{-->}[urr]_{p^l}. & &}
\]
When $l< r$, the degree map $\dr{p^l\colon
P^{2n+1}(p^r)}{}{P^{2n+1}(p^r)}$ is not null homotopic and therefore
neither is its extension $\dr{p^l\colon
S^{2n+1}\{p^r\}}{}{S^{2n+1}\{p^r\}}$. That implies that
$H\exp{S^{2n+1}\{p^r\}}\geq p^r$. On the other hand, for $l=r$, the
degree map $\dr{p^r\colon P^{2n+1}(p^r)}{}{P^{2n+1}(p^r)}$ is null
homotopic. Therefore another choice of $\overline{f}$ extending the
degree $p^r$ map is the trivial map. As $\overline{f}$ is unique,
this implies that $p^r\simeq *$ on $S^{2n+1}\{p^r\}$ and
$H\exp(S^{2n+1}\{p^r\})\leq p^r$.
\end{proof}
\subsection{An upper bound of the exponent for $U$}
Assume that $n>m$. Consider the attaching map $\Sigma\Theta\colon
S^{2n}\lra S^{2m}$ that defines the two-cell complex $P^{2n+1}$. It
is a suspension of a certain homotopy class of $S^{2m-1}$ and
therefore has to be of finite order less than $p^m$
(cf.~\cite{CMNe}). Assume that the order of $\Sigma\Theta$ is $p^l$
for some $l\leq m-1$. Then there is a pushout diagram
\[
\xymatrix{S^{2n}\ar@{=}[d]\ar[r]^{\Sigma\Theta} &
S^{2m}\ar[d]^{p^l}\ar[r] &
P^{2n+1}\ar[d]^a\\
S^{2n}\ar[r]^{*} & S^{2m}\ar[r] & S^{2n+1}\vee S^{2m}}
\]
defining the map $a\colon P^{2n+1}\lra S^{2n+1}\vee S^{2m}$. We
construct a pushout map $\lambda\colon S^{2n+1}\vee S^{2m}\lra
P^{2n+1}$ via the inclusion of the bottom cell $S^{2m}\lra P^{2n+1}$
and the degree $p^l$ map $\dr{p^l\colon P^{2n+1}}{}{P^{2n+1}}$ and
applying the universal property of pushouts. The resulting map
$\lambda$ then gives the commutative diagram
\[
\xymatrix{ S^{2n+1}\vee S^{2m}\ar@{-->}[dr]^{\lambda} &
P^{2n+1}\ar[l]_-a \ar[d]^{p^l}\\
S^{2m}\ar[u]\ar[r] & P^{2n+1}.}
\]
Using the universality of $U$ and the fact that $S^{2n+1}\times
\Omega S^{2m+1}$ is a homotopy associative, homotopy commutative
$H$--space, extend the composite $\ddr{P^{2n+1}}{a}{S^{2n+1}\vee
S^{2m}}{1 \vee E}{S^{2n+1}\times \Omega S^{2m+1}}$ to a unique
$H$-map $\dr{\overline{a}\colon U}{}{S^{2n+1}\times \Omega
S^{2m+1}}.$

As $S^{2n+1}\times \Omega S^{2m+1}$ is a universal space of
$S^{2n+1}\vee S^{2m}$ (cf.~\cite{JGu}), there is a unique
$H$--extension of the composite $\ddr{S^{2n+1}\vee
S^{2m}}{\lambda}{P^{2n+1}}{}{U}$ which we denote by
$\overline{\lambda}\colon S^{2n+1}\times \Omega S^{2m+1}\lra U$.

Composing $\overline{a}$ and $\overline{\lambda}$ together, we
obtain the following diagram
\[
\xymatrix{ P^{2n+1}\ar[rr]^-{p^l}\ar[dr]\ar[dd] & & P^{2n+1}\ar[dd]\\
& S^{2n+1}\vee S^{2m} \ar[ur]\ar[d] &\\
U\ar[r]^-{\overline{a}} & S^{2n+1}\times \Omega S^{2m+1}
\ar[r]^-{\overline{\lambda}} & U.}
\]
By the universality of $U$, the $p^l$ degree map on $P^{2n+1}$
extends to the $p^l$--power map on $U$ that is an $H$--map since $U$
is homotopy commutative. Hence $\overline{\lambda}\overline{a}\simeq
p^l$ giving the factorisation
\[
\xymatrix { U\ar[d]_{\overline{a}}\ar[r]^-{p^l} & U\\
S^{2n+1}\times \Omega S^{2m+1} \ar[ur]_-{\overline{\lambda}}.}
\]
Note that $\exp( S^{2n+1}\times \Omega S^{2m+1})=\max\{\exp
(S^{2n+1}),\exp (S^{2m+1})\}=\max\{p^n,p^m\}=p^n$.

Now it follows that the $p^{l+n}$--power map on $U$ is homotopic to
$\overline{\lambda}\circ p^n\circ \overline{a}$ and therefore it is
trivial. This proves the following proposition.
\begin{proposition}
\label{upper} Let $p^l$ be the order of the attaching map
$\Sigma\Theta$ defining the two-cell complex $P^{2n+1}$. Then
\[H\exp(U)\leq p^{l+n}.\]
\end{proposition}

\subsection{A lower bound of the $H$--homotopy exponent for $U$}
One of the main input data for finding a lower bound of $H\exp(U)$
is the presence of an integer $k\leq m-1$ such that
$p^k\Sigma^2\Theta \simeq *$. That implies the existence of a lift
$h\colon S^{2n+1}\lra U$ of the degree $p^k$ map on $S^{2n+1}$, that
is, a factorisation of the degree $p^k$ map on $S^{2n+1}$ through
$U$
\begin{equation}
\label{MVpullback}\xymatrix{  &
S^{2n+1}\ar[d]^{p^k}\ar[dl]_h\\
 U\ar[r] & S^{2n+1}.}
\end{equation}
Knowing that $\exp (S^{2n+1})=p^n$ (cf. \cite{CMNe}), we can
construct a homotopy class of $U$ so that it is not annihilated by
$p^{n-k-1}$. Therefore we prowed the following proposition.
\begin{proposition}
\label{lower} $p^{n-k}\leq \exp(U)$.
\end{proposition}

Propositions \ref{upper} and \ref{lower} together prove Theorem
\ref{lubounds}.

\end{document}